\documentclass[12pt, oneside]{article}
\usepackage{amsmath,amsfonts,amssymb,amscd,theorem}
\date{}

\newtheorem{thm}{Theorem}
\newtheorem{prop}{Proposition}
\newtheorem{lem}{Lemma}

\begin{document}

\title{Heat flow for horizontal harmonic maps into a class of Carnot-Caratheodory spaces }

\author{J\"urgen Jost and Yi-Hu Yang\thanks{The second named author supported partially
by NSF of China (No. 10171077)}}

\maketitle

\section{Introduction}
Let $X$ and $B$ be two Riemannian manifolds with $\pi: X\to B$ being a Riemannian submersion.
Let $\cal{H}$ be the corresponding horizontal distribution, which is perpendicular to 
the tangent bundle of the fibres of $\pi: X\to B$. Then $X$ (just considered 
as a differentiable manifold),  
together with the distribution $\cal{H}$, forms a so-called {\it Carnot-Caratheodory space} \cite{be}, 
when the Riemannian metric of $X$ is restricted to $\cal{H}$. 
On $X$, as a Carnot-Caratheodory space, can then be defined the notions of Carnot-Caratheodory distance 
(sometimes called sub-Riemannian distance), (minimizing) geodesic, completeness 
(under the Carnot-Caratheodory distance), etc; 
a geodesic is actually a horizontal curve which locally realizes the Carnot-Caratheodory 
distance.
In this note, we always assume that $X$ is complete, as both a Riemannian manifold and
a Carnot-Caratheodory space,  and
the Riemannian submersion $\pi: X\to B$ together with
its horizontal distribution $\cal{H}$ satisfies the following conditions
 
1) the {\it Chow condition}: the vector fields of $\cal{H}$ 
$X_1, X_2, \cdots,$ and their iterated Lie brackets $[X_i, X_j], [[X_i, X_j], X_k], \cdots$ 
span the tangent space $T_xX$ at every point of $X$; 

2) the sectional curvature 
of $X$ (as a Riemannian manifold) in the direction of $\cal{H}$ is non-positive. 

\vskip .3cm
\noindent
{\bf Remark.}
1) The Chow condition guarantees that one has the so-called Hopf-Rinow theorem (cf. \cite{be}): if $X$ is complete under 
the Carnot-Caratheodory metric, then any two points can be joined by a minimizing 
geodesic (under the Carnot-Caratheodory distance); moreover, in any given homotopic class of horizontal
curves connecting two points, there exists a minimzing geodesic 
(under the Carnot-Caratheodory distance) connecting these two points. 2) The Riemannian length
of a horizontal curve is just equal to the Carnot-Caratheodory length by the definitions.
\vskip .3cm
Our interest in this note is to 
study {\it horizontal} maps from a compact Riemannian manifold $M$ into $X$, 
i.e. the image of the derivative of such a map lies in $\cal{H}$.  
We wish to find some 
such maps which furthermore satisfy some differential equation, e.g. harmonic map 
equation, as $X$ is considered as a Riemannian manifold. First of all, let us consider the space of 
smooth maps from $M$ into $X$ which are horizontal and can be connected horizontally to 
a fixed horizontal map $g$, denoted by $B^{\circ}_{g,\cal{H}}(M; X)$; 
it is easy to see that, under a certain suitable metric (defined by using some suitable Sobolev's norm), 
$B^{\circ}_{g,\cal{H}}(M; X)$ can be completed into a Banach manifold, denoted by 
$B_{g,\cal{H}}(M; X)$, which is obviously an infinite dimensional smooth manifold; clearly, its
tangent vectors are just horizontal vector fields of $X$   
(if necessary,they can be considered as  sections of a certain pull-back bundle).
Similarly, considering the space of all maps from 
$M$ into $X$, which are not necessarily horizontal, one can get another Banach manifold, denoted 
by $B(M; X)$, and $B_{g,\cal H}(M; X)$ can be considered as a submanifold of $B(M; X)$. It should be pointed 
out that these Banach manifolds may not be connected 
(but clearly are locally connected), this does not however affect our following discussion.
Let $\cal X$ be a vector field of $B(M; X)$ along $B_{g,\cal H}(M; X)$. Corresponding to the horizontal
distribution $\cal H$, one has an orthogonal projection to $\cal H$, still denoted by $\cal H$. 
Accordingly, one can also
define the projection of $\cal X$, denoted by $\cal{H}\cal{X}$, which is a vector field of $B_{g,\cal H}(M; X)$
and the value of which at any point of $B_{g,\cal{H}}(M; X)$ is actually  
a horizontal vector field of $X$ 
(again, if necessary, it can be considered as a section of a certain pull-back bundle).

In this note, we first give some examples of Carnot-Caratheodory spaces, in which we are really interested.  
These spaces are actually a class of (locally) complex homogeneous manifolds which fibre over the corresponding
symmetric spaces of noncompact type and the fiberations are Riemannian submersion under the standard invariant metrics.
We will show that this class of spaces satisfies the Conditions 1) and 2) above; 
on the other hand, such homogeneous spaces, as Riemannian manifolds, are complete 
and by the definition of Carnot-Caratheodory distance, 
the Riemannian distance is not greater than the Carnot-Caratheodory distance, 
so this class of homogeneous spaces, as Carnot-Caratheodory spaces, are also complete under the corresponding 
Carnot-Caratheodory distance. 
Thus we can apply the Banach spaces defined 
above to this class  of homogeneous complex manifolds.
We next consider the following heat flow from $M\times [0, \infty)$ into $X$
\[
(*)~~~~~~~~~~~~~~~~~~~~~~~~~~~~~~~~~~~~~~{\cal H}\tau(u) - {\frac{\partial{u}}{\partial t}} = 0,~~~~
~~~~~~~~~~~~~~~~~~~~~~~~~~~~~~~~~~~~~~~~~~~~
\]
with the initial data $u(\cdot, 0)=g(\cdot)$, here $\tau(u)$ is the stress-energy tensor of $u$ with respect to
the space variable, $g(\cdot)$ is a smooth horizontal map. We show that one can always deform
horizontally any smooth horizontal map into a horizontal harmonic map. It is worth noting that 
the operator ${\cal H}\tau$, as applied to the Banach space $B(M; X)$, is not elliptic in general, 
but if applied to the Banach space $B_{g,{\cal H}}(M; X)$, it is indeed elliptic, 
i.e., the symbol of its linearization is an isomorphism from the horizontal tangent subbundle of $X$ to itself, 
and hence one can apply the implicit
function theorem to the Banach space $B_{g,{\cal H}}(M; X)$ to obtain
the short-time existence of a (unique) solution of $(*)$
with the initial map $g$.

\vskip .5cm
\noindent
{\bf Acknowlegements}. The second named author wants to thank Dr Guofang Wang for valuable discussion.
The main part of this work was completed when he was visiting 
the Max-Planck-Institute for Mathematics in the Sciences. 
He would like also to thank the institute for its hospitality and good working conditions.

\section{A class of Carnot-Caratheodory spaces}
In this section, we will show some concrete examples for Carnot-Caratheodory spaces, which 
are actually the objects in which we are really interested.
These examples are a class of (locally) complex homogeneous manifolds \cite{gsch, jy}: 
Let $G$ be a connected noncompact real semisimple Lie group satisfying
that it has a compact Cartan subgroup; 
as a consequence, if $K$ is a maximal compact subgroup of $G$, then $G$ and $K$ have the same rank;
moreover $G/K$ is not a Hermitian symmetric space. Denote such a Cartan 
subgroup by $H$, and choose a suitable subgroup $Z$ of $K$ containing $H$, which is actually the centralizer 
in $G$ of a certain circle subgroup $T$ of $H$. Taking the quotients $G/Z$ and $G/K$, one has then that
$G/Z$ is a homogeneous complex manifold and $G/K$ is a symmetric space of noncompact type; 
moreover $G/Z$ is a fiberation
over $G/K$ with the fiber $K/Z$; under the standard invariant metrics \cite{jy}, the fibration 
$\pi: G/K\to G/Z$ is a Riemannian submersion, and hence it has a horizontal distribution $\cal H$, 
which satisfies all the assumptions mentioned in the preceeding section, as shown in the following. 
Let $\Gamma$ be a discrete subgroup of $G$. Because of the discreteness of $\Gamma$ and the compactness of $K$, 
one can assume that $\Gamma\cap K=\emptyset$. Thus we have the Riemannian submersion
$\Gamma\setminus G/Z\to\Gamma\setminus G/K $. Similarly, one has the horizontal distribution 
which is the discrete quotient of $\cal H$ and hence also satisfies the assumption in the preceeding section, 
denoted by $\cal H'$. 
In the remaining part of this section, we will show that the distribution $\cal H$, and hence $\cal H'$,
does satisfy those assumptions.
First, we check the assumption for sectional curvature in 
the horizontal direction $\cal H$; actually, one generally has the following
\begin{prop}
Let $\pi: X\to B$ be a Riemannian submersion. If $B$ has non-positive sectional curvature, then $X$, 
in the horizontal direction $\cal{H}$, also has non-positive sectional curvature.   
\end{prop}

Since $G/K$ is a symmetric space of noncompact type, so it, and hence $G/Z$ 
in the horizontal direction $\cal H$, has non-positive sectional curvature.

\vskip .3cm
\noindent
{\bf Proof of Proposition 1.} The proof is a simple consequence of the O'Neill formulae: Denote the curvature tensors of 
$X$ and $B$ by $R$ and $R'$ respectively; then one of O'Neill's formulae says, for horizontal tangent 
vectors $Y, Z, U, V$ of $X$,
\begin{eqnarray*}
&&<R(Y, Z)U, V> = <R'(Y, Z)U, V> - 2<A(Y, Z), A(U, V)> \\
&&+ <A(Z, U), A(Y, V)> -<A(Y, U), A(Z, V)>.
\end{eqnarray*}
Here, $Y, \cdots$ are also regarded as tangent vectors of $B$; the definition of $A$ refers to 
the proof of the Lemma 2 in the next section; the key point is that $A$ is 
skew-symmetric with respect to horizontal vectors. So
\begin{eqnarray*}
&&<R(Y, Z)Y, Z> = <R'(Y, Z)Y, Z> - 2<A(Y, Z), A(Y, Z)> + \\
&&<A(Z, Y), A(Y, Z)> = <R'(Y, Z)Y, Z> - 3<A(Y, Z), A(Y, Z)> \le 0.
\end{eqnarray*}

We now turn to check the Chow condition.
By Cartan's classification theorem for simple groups \cite{hel}, 
the simple Lie groups satisfying the conditions stated in the beginning of this section are as follows :
\begin{eqnarray*}
&~~~~~~~~~~~~~~~~{\text {SO}}(p, 2q)~q\ge 2~&{\bf e}_{8(8)}~~~~~~~~~~~~~~~~~~~~~~~~~~~~~~~~~~~~~~~~~~~~~\\
&~~~~{\text Sp}(p, q)&{\bf e}_{8(-24)}~~~~~~~~\\
&{\bf e}_{6(2)}&{\bf f}_{4(4)}~~~~~~~~\\
&{\bf e}_{7(7)}&{\bf f}_{4(-20)}~~~~\\
&~~{\bf e}_{7(-5)}&{\bf g}_{2(2)}
\end{eqnarray*}
The above list is called  {\it groups of Hodge type but not of Hermitian type} in Simpson's paper\cite{si}. 
In order to show the Chow condition, we can actually turn the problem into a Lie-theoretic problem.
To this end, we first need to give the relation between the Lie bracket of left invariant vector fields 
and the Lie bracket of the Lie algebras in question when considering left invariant vector fields 
as elements of the Lie algebra. We use the notations of \cite{kn}. Denote the Lie algebra of $G$ and $Z$ by 
$\mathfrak{g}$ and $\mathfrak{z}$ respectively, then it is easy to see that we have a direct 
sum decomposition of vector spaces
\[
\mathfrak{g} = \mathfrak{z} + \mathfrak{m}
\]
with $[\mathfrak{z}, \mathfrak{m}]\subset\mathfrak{m}$. Here $\mathfrak{m}$ can be identified with
the tangent space of $G/Z$ at the origin or the set of all $G$-invariant vector fields on $G/Z$. 
Theorem 2.10 of \cite{kn} tells us that there exists a unique torsion-free $G$-invariant affine 
connection $\nabla$ with
\[
\nabla_YZ = {\frac 1 2}[Y, Z]_{\mathfrak{m}}, ~{\text{for}}~Y, Z\in {\mathfrak{m}},
\]
here by $Y, Z$ on the left-hand side we mean vector fields on $G/Z$ while $Y, Z$ on the right-hand 
side mean elements in $\mathfrak{g}$; $[Y, Z]_{\mathfrak{m}}$ denotes the $\mathfrak{m}$-component of 
$[Y, Z]$. Thus, one has
\[
[Y, Z] = [Y, Z]_{\mathfrak{m}},
\]
here by the left-hand side we mean the Lie bracket of vector fields; afterwards we will not point out
this since it should be clear from the context. As before, one has a Cartan subgroup $H$ contained in
$Z$, the Lie algebra of which is a maximal abelian subalgebra, denoted by $\mathfrak{h}$. 
Consider the Cartan decomposition 
$\mathfrak{g}=\mathfrak{k}+\mathfrak{p}$, here $\mathfrak{k}$ is the Lie algebra of $K$. 
We then have the following relations $\mathfrak{h}\subset\mathfrak{z}\subset\mathfrak{k}\subset\mathfrak{g}$
and $\mathfrak{p}\subset\mathfrak{m}$. Again, $\mathfrak{p}$, as a vector subspace of $\mathfrak{m}$, 
can be indentified with the horizontal tangent subspace at the origin with respect to the Riemannian 
submersion $G/Z\to G/K$ and its left translation forms the horizontal distribution $\cal H$ of the Riemannian submersion; furthermore, its elements can be identified with $G$-invariant horizontal
vector fields of $G/K$. By the previous relation of two Lie brackets, in order to show that 
the horizontal distribution $\cal H$ satisfies the Chow condition, 
it is sufficient to show 
that $\mathfrak{p}$ and $[\mathfrak{p}, \mathfrak{p}]$ span $\mathfrak{m}$. 
To this end, we use the root system of the complexification $\mathfrak{g}^{\bf C}$ of $\mathfrak{g}$ 
correspoding to the Cartan subalgebra $\mathfrak{h}$. Let $\Delta$ be the root system of 
$\mathfrak{g}^{\bf C}$ with respect to $\mathfrak{h}$, $\mathfrak{g}^{\alpha}$ the root 
space corresponding to $\alpha\in\Delta$, 
$\mathfrak{g}^{\bf C}=\mathfrak{k}^{\bf C}+\mathfrak{p}^{\bf C}$ the Cartan decomposition, $\theta$
the Cartan involution, $\sigma$ the conjugation of $\mathfrak{g}^{\bf C}$ with respect to $\mathfrak{g}$.
Since $\mathfrak{h}$ lies in $\mathfrak{k}$ while $[\mathfrak{k}, \mathfrak{k}]\subset\mathfrak{k}$ and
$[\mathfrak{k}, \mathfrak{p}]\subset\mathfrak{p}$, so the root space $\mathfrak{g}^{\alpha}$ lies in either
$\mathfrak{k}^{\bf C}$ or $\mathfrak{p}^{\bf C}$. In the first case,
we call $\alpha$ a {\it compact root}; 
denote the set of all compact roots 
by $\Delta(\mathfrak{k})$; in the 
last case, a {\it noncompact root}; denote the set of noncompact roots by $\Delta(\mathfrak{p})$.
On the other hand, we also have the direct sum decomposition for vector spaces 
$\mathfrak{g}=\mathfrak{h}+\mathfrak{m}'$, obviously $\mathfrak{m}\subset\mathfrak{m}'$; furthermore
one has the direct sum $\mathfrak{m}'=\mathfrak{k}'+\mathfrak{p}$ with 
$\mathfrak{h}+\mathfrak{k}'=\mathfrak{k}$. So if we can show that 
$[\mathfrak{p}, \mathfrak{p}]=\mathfrak{k}'$, equivalently 
$[\mathfrak{p}^{\bf C}, \mathfrak{p}^{\bf C}]=\mathfrak{k'}^{\bf C}$, 
then the Chow condition is obtained.
>From the root theory, we has
\[
\mathfrak{k'}^{\bf C}=\sum_{\alpha\in\Delta(\mathfrak{k})}{\mathfrak{g}}^{\alpha} ~~{\text{and}}~~
\mathfrak{p}^{\bf C}=\sum_{\alpha\in\Delta(\mathfrak{p})}{\mathfrak{g}}^{\alpha}.
\]
Note that $\sigma({\mathfrak{g}}^{\alpha})={\mathfrak{g}}^{-\alpha}$ 
while $\sigma(\mathfrak{k'}^{\bf C})=\mathfrak{k'}^{\bf C}$ and
$\sigma(\mathfrak{p}^{\bf C})=\mathfrak{p}^{\bf C}$,
so if $\alpha\in\Delta(\mathfrak{k})$ (resp. $\Delta(\mathfrak{p})$), then so is $-\alpha$.
We now state the following 
\begin{prop} 
For any root $\alpha\in\Delta(\mathfrak{k})$, there exist two noncompact roots
$\beta$ and $\gamma$ with $\beta+\gamma=\alpha$. 
\end{prop}
Clearly if the proposition is true,
then the Chow codition is obtained.
In the following, we will case by case write down compact roots 
and noncompact roots of $\mathfrak{g}^{\bf C}$ for the above simple groups list and then easily 
check that the above assertion is true.

\vskip .5cm
\noindent
$SO(p, 2q), q\ge 2$: we have two cases to consider. $SO(2p, 2q), p, q\ge 2$: It is the noncompact 
real form of $SO(2(p+q), {\bf C})$ with the maximal compact subgroup $K=SO(2p)\times SO(2q)$.
The root system of $\mathfrak{so}(2(p+q), {\bf C})$ is 
$D_{p+q}=\{\pm e_i\pm e_j, 1\le i<j\le p+q\}$, here $\{e_i\}$
is the standard basis of ${\bf R}^{p+q}$, while the root systems of $\mathfrak{so}(2p, {\bf C})$ and 
$\mathfrak{so}(2q, {\bf C})$, embedded in $D_{p+q}$, are 
\[
D_{p}=\{\pm e_i\pm e_j, 1\le i<j\le p\}
\] 
and
\[
D_{q}=\{\pm e_i\pm e_j, p+1\le i<j\le p+q\}
\] 
respectively. 
Therefore, corresponding to the noncompact 
real form $SO(2p, 2q)$ and its compact Cartan subalgebra, $\mathfrak{so}(2(p+q), {\bf C})$ has noncompact roots
\[
\{\pm e_i\pm e_j, 1\le i\le p, p+1\le j\le p+q\};
\] 
the second case is $SO(2p+1, 2q), p, q\ge 2$: 
it is the noncompact real form of $SO(2(p+q)+1, {\bf C})$ with the maximal compact 
subgroup $K=SO(2p+1)\times SO(2q)$. 
The root system of $\mathfrak{so}(2(p+q)+1, {\bf C})$ is 
$B_{p+q}=\{\pm e_i, \pm e_i\pm e_j, 1\le i, j\le p+q, i\neq j\}$
while the root systems of $\mathfrak{so}(2p+1, {\bf C})$ and 
$\mathfrak{so}(2q, {\bf C})$, embedded in $B_{p+q}$, are
\[ 
B_{p}=\{\pm e_i, \pm e_i\pm e_j, 1\le i,j\le p, i\neq j\} 
\]
and
\[
D_{q}=\{\pm e_i\pm e_j, p+1\le i<j\le p+q\} 
\]
respectively. 
Therefore, corresponding to the noncompact 
real form $SO(2p+1, 2q)$ and its compact Cartan subalgebra, $\mathfrak{so}(2(p+q)+1, {\bf C})$ has noncompact roots
\[
\{\pm e_i\pm e_j, \pm e_j, 1\le i\le p, p+1\le j\le p+q\}.
\]

\vskip .5cm
\noindent
$Sp(p, q)$: It is the noncompact real form of $Sp(p+q, {\bf C})$ with the maximal compact subgroup
$K=Sp(p)\times Sp(q)$. The root system of $\mathfrak{sp}(p+q, {\bf C})$ is 
$C_{p+q}=\{\pm 2e_i, \pm e_i\pm e_j, 1\le i, j\le p+q, i\neq j\}$,
while the root systems of $\mathfrak{sp}(p, {\bf C})$ and $\mathfrak{sp}(q, {\bf C})$, embedded in $C_{p+q}$,
are 
\[
C_p=\{\pm 2e_i, \pm e_i\pm e_j, 1\le i, j\le p, i\neq j\} 
\]
and
\[
C_q=\{\pm 2e_i, \pm e_i\pm e_j, p+1\le i, j\le p+q, i\neq j\}
\]
respectively; 
therefore, corresponding to the noncompact 
real form $Sp(p, q)$ and its compact Cartan subalgebra, $\mathfrak{sp}(p+q, {\bf C})$ has noncompact roots
\[
\{\pm e_i\pm e_j, 1\le i\le p, p+1\le j\le p+q\}.
\]

\vskip .5cm
\noindent
${\bf e}_{6(2)}$: It is the noncompact real form of ${\bf e}_6$ with the maximal compact subgroup
$K=SU(6)\times SU(2)$. The root system of ${\bf e}_6$ is 
\begin{eqnarray*}
E_6&=&\{e_i-e_j, i\neq j, 1\le i, j\le 6\}\cup\{\pm(e_7-e_8)\}\cup  \\
&&\{{\frac 1 2}(e_{\sigma(1)}+e_{\sigma(2)}+e_{\sigma(3)}
-e_{\sigma(4)}-e_{\sigma(5)}-e_{\sigma(6)}\pm(e_7-e_8)), \sigma\in P(6)\},
\end{eqnarray*}
where $P(6)$ is the permutation group of $\{1,2,3,4,5,6\}$. 
The root system of $\mathfrak{sl}(6, {\bf C})+\mathfrak{sl}(2, {\bf C})$, embedded in $E_6$, is
\[
A_5+A_1=\{e_i-e_j, i\neq j, 1\le i, j\le 6\}\bigcup\{\pm(e_7-e_8)\}.
\] 
Thus, corresponding to the noncompact 
real form ${\bf e}_{6(2)}$ and its compact Cartan subalgebra, ${\bf e}_6$ has noncompact roots
\[
\{{\frac 1 2}(e_{\sigma(1)}+e_{\sigma(2)}+e_{\sigma(3)}
-e_{\sigma(4)}-e_{\sigma(5)}-e_{\sigma(6)}\pm(e_7-e_8)), \sigma\in P(6)\}.
\]

\vskip .5cm
\noindent
${\bf e}_{7(7)}$: It is the noncompact real form of ${\bf e}_7$ with the maximal compact subgroup
$K=SU(8)$. The root system of ${\bf e}_7$ is 
\begin{eqnarray*}
E_7&=&\{e_i-e_j, 1\le i, j\le 8, i\neq j\}\cup  \\ 
&&\{{\frac 1 2}(e_{\sigma(1)}+e_{\sigma(2)}+e_{\sigma(3)}+e_{\sigma(4)}-
e_{\sigma(5)}-e_{\sigma(6)}-e_{\sigma(7)}-e_{\sigma(8)}), \sigma\in P(8)\}, 
\end{eqnarray*}
here $P(8)$ is the permutation group of $\{1,2,3,4,5,6,7,8\}$. 
The root system of $\mathfrak{sl}(8, {\bf C})$, embedded in $E_7$, is
\[
A_7=\{e_i-e_j, i\neq j, 1\le i, j\le 8\}.
\] 
Thus, corresponding to the noncompact 
real form ${\bf e}_{7(7)}$ and its compact Cartan subalgebra, ${\bf e}_7$ has noncompact roots
\[
\{{\frac 1 2}(e_{\sigma(1)}+e_{\sigma(2)}+e_{\sigma(3)}+e_{\sigma(4)}-
e_{\sigma(5)}-e_{\sigma(6)}-e_{\sigma(7)}-e_{\sigma(8)}), \sigma\in P(8)\}.
\]

\vskip .5cm
\noindent
${\bf e}_{7(-5)}$: It is the noncompact real form of ${\bf e}_7$ with the maximal compact subgroup
$K=SO(12)\times SU(2)$. The root system of ${\bf e}_7$ is 
\begin{eqnarray*}
E_7&=&\{e_i-e_j, 1\le i, j\le 8, i\neq j\}\cup  \\ 
&&\{{\frac 1 2}(e_{\sigma(1)}+e_{\sigma(2)}+e_{\sigma(3)}+e_{\sigma(4)}-
e_{\sigma(5)}-e_{\sigma(6)}-e_{\sigma(7)}-e_{\sigma(8)}), \sigma\in P(8)\}, 
\end{eqnarray*}
here $P(8)$ is the permutation group of $\{1,2,3,4,5,6,7,8\}$. 
The root system of $\mathfrak{so}(12, {\bf C})+\mathfrak{sl}(2, {\bf C})$, embedded in $E_7$, is
\begin{eqnarray*}
D_6+A_1&=&\{e_i-e_j, 1\le i, j\le 6, i\neq j\}\cup   \\
&&\{\pm{\frac 1 2}(e_{\sigma(1)}+e_{\sigma(2)}+e_{\sigma(3)}+e_{\sigma(4)}-
e_{\sigma(5)}-e_{\sigma(6)}-e_7-e_8)\}\bigcup  \\
&&\{\pm(e_7-e_8)\}.
\end{eqnarray*}

\noindent
(Note that if letting ${\bf R}^n$ have the standard basis $\{f_1, \cdots f_n\}$, 
$D_n=\{\pm f_i\pm f_j, i\neq j\}$; so we need to construct an isomorphism between $D_6$ and 
$\{e_i-e_j, 1\le i, j\le 6, i\neq j\}\cup   
\{\pm{\frac 1 2}(e_{\sigma(1)}+e_{\sigma(2)}+e_{\sigma(3)}+e_{\sigma(4)}-
e_{\sigma(5)}-e_{\sigma(6)}-e_7-e_8)\}$. This is done by the uniqueness: 
$\{e_i-e_j, 1\le i, j\le 6, i\neq j\}\cup   
\{\pm{\frac 1 2}(e_{\sigma(1)}+e_{\sigma(2)}+e_{\sigma(3)}+e_{\sigma(4)}-
e_{\sigma(5)}-e_{\sigma(6)}-e_7-e_8)\}$ indeed is a root system of cardinality $60$; on the other hand, the root 
system of cardinality $60$ is only $D_n$ by the Cartan classification theorem.)

\noindent
Therefore, corresponding to the noncompact 
real form ${\bf e}_{7(-5)}$ and its compact Cartan subalgebra, ${\bf e}_7$ has noncompact roots
\begin{eqnarray*}
&&\{\pm(e_i-e_7), \pm(e_i-e_8), 1\le i\le 6\}\cup  \\
&&\{\pm{\frac 1 2}(e_{\sigma(1)}+e_{\sigma(2)}+e_{\sigma(3)}-e_{\sigma(4)}-
e_{\sigma(5)}-e_{\sigma(6)}+e_7-e_8)\}.
\end{eqnarray*}

\vskip .5cm
\noindent
${\bf e}_{8(8)}$: It is the noncompact real form of ${\bf e}_8$ with the maximal compact subgroup
$K=SO(16)$. The root system of ${\bf e}_8$ is 
\[
E_8=\{\pm e_i\pm e_j, {\frac 1 2}\sum_{i=1}^{8}(-1)^{m(i)}e_i ~{\text{with}}~ 
\sum m(i) ~{\text{being even}}, 1\le i, j\le 8\},
\]
where $m(i)$ is $0$ or $1$. The root system of ${\mathfrak{so}}(16, {\bf C})$, embedded in $E_8$, is 
$D_8=\{\pm e_i\pm e_j, 1\le i, j\le 8\}$. Therefore, corresponding to the noncompact 
real form ${\bf e}_{8(8)}$ and its compact Cartan subalgebra, ${\bf e}_8$ has noncompact roots
\[
\{{\frac 1 2}\sum_{i=1}^{8}(-1)^{m(i)}e_i ~{\text{with}}~ 
\sum m(i) ~{\text{being even}}, 1\le i, j\le 8\}.
\]

\vskip .5cm
\noindent
${\bf e}_{8(-24)}$: It is the noncompact real form of ${\bf e}_8$ with the maximal compact subgroup
$K={\bf e}_{7(-133)}\times SU(2)$. The root system of ${\bf e}_8$ is 
\[
E_8=\{\pm e_i\pm e_j, {\frac 1 2}\sum_{i=1}^{8}(-1)^{m(i)}e_i ~{\text{with}}~ 
\sum m(i) ~{\text{being even}}, 1\le i, j\le 8\}. 
\]
The root system of ${\bf e}_7+sl(2, {\bf C})$, embedded in $E_8$, is 
\begin{eqnarray*}
&&E_7+A_1=\{e_i-e_j, 1\le i, j\le 8, i\neq j\}\cup\\ 
&&\{{\frac 1 2}(e_{\sigma(1)}+e_{\sigma(2)}+e_{\sigma(3)}+e_{\sigma(4)}-
e_{\sigma(5)}-e_{\sigma(6)}-e_{\sigma(7)}-e_{\sigma(8)}), \sigma\in P(8)\}\bigcup  \\
&&\{\pm{\frac 1 2}(e_1+e_2+ \cdots +e_8)\}. 
\end{eqnarray*}
Thus, corresponding to the noncompact 
real form ${\bf e}_{8(-24)}$ and its compact Cartan subalgebra, ${\bf e}_8$ has noncompact roots
\begin{eqnarray*}
&&\{\pm(e_i+e_j), 1\le i<j\le 8\}\cup\\
&&\{{\frac 1 2}(e_{\sigma(1)}+e_{\sigma(2)}+e_{\sigma(3)}+e_{\sigma(4)}+
e_{\sigma(5)}+e_{\sigma(6)}-e_{\sigma(7)}-e_{\sigma(8)}), \sigma\in P(8)\}.
\end{eqnarray*}

\vskip .5cm
\noindent
${\bf f}_{4(4)}$: It is the noncompact real form of ${\bf f}_4$ with the maximal compact subgroup
$Sp(3)\times SU(2)$. The root system of ${\bf f}_4$ is
\[
F_4=\{\pm e_i, \pm e_i\pm e_j~ (1\le i, j\le 4, i\neq j), {\frac 1 2}(\pm e_1\pm e_2\pm e_3\pm e_4)\};
\]
while the root system of $\mathfrak{sp}(3, {\bf C})+\mathfrak{sl}(2, {\bf C})$, 
embedded in $F_4$, is
\[
C_3+A_1=\{\pm 2f_i, \pm f_i\pm f_j, 1\le i, j\le 3, i\neq j\}\bigcup\{\pm(e_3+e_4)\}
\]
where $f_1={\frac 1 2}(e_1-e_2),f_2={\frac 1 2}(e_1+e_2), f_3={\frac 1 2}(e_3-e_4)$.
Thus, corresponding to the noncompact 
real form ${\bf f}_{4(4)}$ and its compact Cartan subalgebra, ${\bf f}_4$ has noncompact roots
\begin{eqnarray*}
&&\{\pm e_3, \pm e_4, \pm e_i\pm e_j, i=1, 2,j=3, 4\}\cup  \\
&&\{{\frac 1 2}(\pm(e_1-e_2)\pm(e_3+e_4)), {\frac 1 2}(\pm(e_1+e_2)\pm(e_3+e_4))\}.
\end{eqnarray*}

\vskip .5cm
\noindent
${\bf f}_{4(-20)}$: It is the noncompact real form of ${\bf f}_4$ with the maximal compact subgroup
$SO(9)$. The root system of ${\bf f}_4$ is
\[
F_4=\{\pm e_i, \pm e_i\pm e_j~ (1\le i, j\le 4, i\neq j), {\frac 1 2}(\pm e_1\pm e_2\pm e_3\pm e_4)\};
\]
while the root system of $\mathfrak{so}(9, {\bf C})$, embedded in $F_4$, is 
$B_4=\{\pm e_i, \pm e_i\pm e_j, 1\le i, j\le 4, i\neq j\}$.
Therefore, corresponding to the noncompact 
real form ${\bf f}_{4(-20)}$ and its compact Cartan subalgebra, ${\bf f}_4$ has noncompact roots
\[
\{{\frac 1 2}(\pm e_1\pm e_2\pm e_3\pm e_4)\}.
\]

\vskip .5cm
\noindent
${\bf g}_{2(2)}$: It is the noncompact real form of ${\bf g}_2$ with the maximal compact subgroup
$SU(2)\times SU(2)$. The root system of ${\bf g}_2$ is 
\[
G_2=\{\pm\alpha, \pm\beta, \pm(\alpha+\beta), \pm(2\alpha+\beta), \pm(3\alpha+\beta), \pm(3\alpha+2\beta)\},
\]
where $\alpha=e_1, \beta=-{\frac 3 2}e_1+{\frac{\sqrt{3}}2}e_2$; while the root system of 
${\mathfrak{sl}}(2, {\bf C})+{\mathfrak{sl}}(2, {\bf C})$, embedded in $G_2$,
is $A_1+A_1=\{\pm\beta\}\bigcup\{\pm(2\alpha+\beta)\}$. Therefore the noncompact root system is
\[
\{\pm\alpha, \pm(\alpha+\beta), \pm(3\alpha+\beta), \pm(3\alpha+2\beta)\}.
\]
Summing the above all up, it is easy to check that the noncompact roots can generate the
compact roots, i.e. for any compact root 
$\alpha$ there exist two noncompact roots $\beta$ and $\gamma$ satisfying $\alpha=\beta+\gamma$.

\section{Heat flow for horizontal harmonic maps}
Let $\pi: X\to B$ be a Riemannian submersion, $\cal H$ the corresponding horizontal distribution, 
and $M$ a compact Riemannian manifold. Assume that $\pi: X\to B$ satisfies the conditions stated in the Introduction, 
i.e. the Chow condition and $B$ having non-positive sectional curvature 
and that $X$ is complete under both the Carnot-Caratheodory distance and the Riemannian metric.  
Consider the following heat equation on $M$
\[
(*)~~~~~~~~~~~~~~~~~~~~~~~~~~~~~~~~~~~~~~{\cal H}\tau(u) - {\frac{\partial{u}}{\partial t}} = 0,~~~~
~~~~~~~~~~~~~~~~~~~~~~~~~~~~~~~~~~~~~~~~~~~~
\]
where $\cal H$ represents the projection to $\cal H$, and $\tau$ is
the tension field (nonlinear Laplacian) of $u$.
Assume that $u$ has initial data $u(\cdot, 0) = g(\cdot)$. We always
assume that $g$ is a smooth horizontal 
map from $M$ to $X$. We wish to obtain some horizontal harmonic map from 
$M$ into $X$ by solving the above heat equation for the initial data $g$, 
when $X$ is considered as a Riemannian manifold. 
\begin{lem}
There exists a positive number $T$, such that the equation $(*)$ with the initial data $g$ has a smooth solution 
$u(x, t)$ for $t\in [0, T)$ satisfying $u(\cdot, t)\in B_{g,\cal H}(M; X)$. 
Furthermore, if $u(x, t)$ is a solution of $(*)$ with $u(\cdot, 0)=g(\cdot)$ 
for $t\in [0, T'), T'>0$, then $u(\cdot, t)\in B_{g, \cal{H}}(M; X)$ and hence 
${\frac{u(\cdot, t)}{\partial t}}$ is a horizontal tangent vector field of $X$ for any $t\in [0, T')$. 
\end{lem}    
\noindent
{\bf Proof.} The first part of the lemma is essentially a standard result if one restricts the problem
to the space $B_{g,{\cal H}}(M; X)$: The symbol of the linearization of the operator ${\cal H}\tau$ is just
an isomorphism from the horizontal tangent subbundle of $X$ to itself,
so ${\cal H}\tau$ is elliptic. Thus one can still apply the implicit function theorem to the present case,
as one applies the implicit function theorem to the usual harmonic map heat flow, to obtain the short-time
existence.
As for the second part, it is also easy to see from
the following discussion.
Since ${\cal H}\tau(u)$ is a horizontal vector on $B(M; X)$, i.e. a horizontal vector field on $X$, so 
${\frac{\partial u}{\partial t}}$ is also horizontal. Fix a point $x\in M$ and take arbitrarily
a curve $\gamma(s)$ starting from $x$ for $s\in [0, s_0]$ and a vertical tangent vector $V$ at $g(x)$, translate
parallelly $V$ along the $t$-curve $u(x, t)$ and then the $s$-curves $u(\gamma(s), t)$, still denoted by $V$.
Note that $V$ is not necessarily parallel, even not continuous, along the $t$-curves $u(\gamma(s), t)$ for $s\neq 0$.
Compute 
${\frac{\partial}{\partial t}}<{\frac{\partial}{\partial s}}u(\gamma(0), t), V>$
\begin{eqnarray*}
&&{\frac{\partial}{\partial t}}<{\frac{\partial}{\partial s}}u(\gamma(0), t), V>
=<\nabla_{{\frac{\partial}{\partial t}}}{\frac{\partial}{\partial s}}u(\gamma(0), t), V>\\
&=&<\nabla_{{\frac{\partial}{\partial s}}}{\frac{\partial}{\partial t}}u(\gamma(0), t), V>
={\frac{\partial}{\partial s}}<{\frac{\partial}{\partial t}}u(\gamma(0), t), V>=0.
\end{eqnarray*}
Since $<{\frac{\partial}{\partial s}}u(\gamma(0), t), V>|_{t=0}=<{\frac{\partial}{\partial s}}g(\gamma(0)), V>=0$,
so $<{\frac{\partial}{\partial s}}u(\gamma(0), t), V>=0$. Thus $u(\cdot, t)$ is horizontal. 
Then, the horizontality of ${\frac{u(\cdot, t)}{\partial t}}$ implies $u(\cdot, t)\in B_{g, \cal{H}}(M; X)$.
The lemma is obtained.

\vskip .3cm
Let $e(u)(x, t)={\frac 1 2}|\nabla u|^2(x, t)$ be the energy density of $u(\cdot, t)$ for $t\in [0, T)$.
Denote the Laplace operator of $M$ by $\Delta$ and take $\{e_i\}$ as a normal frame of $M$;
denote the Ricci tensor of $M$ by ${\text{Ric}}^M$ and the curvature tensor of $X$ by $R^X$.
By $\cal V$ we mean to take the vertical component of vectors.
Then compute $(\Delta - {\frac{\partial}{\partial t}})e(u)$:
\begin{eqnarray*}
(\Delta - {\frac{\partial}{\partial t}})e(u)
&=& <\nabla_{e_i}\nabla_{e_i}{\text d}u, {\text d}u> + |\nabla{\text d}u|^2 - 
<\nabla{\frac{\partial u}{\partial t}}, {\text d}u>\\
&=& <\nabla({\cal V}\tau(u)), {\text d}u> + |\nabla{\text d}u|^2 
+ <{\text{Ric}}^M({\text d}u(e_i), {\text d}u(e_i))>  \\
&&- <R^X({\text d}u(e_i), {\text d}u(e_j)){\text d}u(e_i), {\text d}u(e_j)>\\
&=& - |{\cal V}\tau(u)|^2 + |\nabla{\text d}u|^2 
+ <{\text{Ric}}^M({\text d}u(e_i), {\text d}u(e_i))>  \\
&&- <R^X({\text d}u(e_i), {\text d}u(e_j)){\text d}u(e_i), {\text d}u(e_j)>.
\end{eqnarray*}
In the second equality above we used the Weitzenb\"ock formula and the equation $(*)$;
in the last equality we used the horizontality of $u$.
The following observation is important for the present study.

\begin{lem}
Let $\pi: X\to B$ be a Riemannian submersion. Then, for any horizontal map
$u$ from a Riemannian manifold $M$ into $X$, the vertical part ${\cal V}\tau(u)$ of its stress-energy tensor
$\tau(u)$ vanishes. 
\end{lem}

\vskip .3cm
\noindent
{\bf Remark.} Since the horizontal distribution $\cal H$ is generally not integrable, so the vertical part
of the Hessian of a horizontal map $u$ does not necessarily vanish.

\vskip .3cm
\noindent
{\bf Proof.} We first review an idea of B. O'Neill \cite{bes, on}. According to 
B. O'Neill, one can define a type $(2, 1)$-tensor
field on $X$, denoted by $A$, as follows: for any two vector fields $Y, Z$ on $X$,
\[
A(Y, Z) = {\cal H}\nabla_{{\cal H}Y}{\cal V}Z + {\cal V}\nabla_{{\cal H}Y}{\cal H}Z,
\]
here $\cal H$ and $\cal V$ mean taking the horizontal part and the vertical part 
respectively, as mentioned before. An easy calculation shows that $A$ indeed is a tensor 
field on $X$, namely, the value of $A(Y, Z)$ at any fixed point $x$
depends only on the values of
$Y$ and $Z$ at $x$, although its definition does depend on the value of $Y$ and $Z$ on a small 
neighborhood of $x$; moreover, it has the following key property 
(here we state slightly  more than we actually need):
\[
A(Y, Z) = -A(Z, Y) = {\frac 1 2}{\cal V}[Y, Z]
\]
for any two horizontal vectors $Y$ and $Z$. The proof of this property is simple: It is sufficient to show
$A(Y, Y)=0$. Namely if this is the case, $A(Y+Z, Y+Z)=A(Y, Z)+A(Z, Y)=0$; and, by the definition of $A$,
\[
{\cal V}[Y, Z] = {\cal V}\nabla_YZ - {\cal V}\nabla_ZY = A(Y, Z) - A(Z, Y).
\] 
Since $A$ is a tensor,  one can take the horizontal vector field $Y$ being 
the unique lift of a vector field $Y'$ on $B$, i.e. $\pi_*(Y)=Y'$. Let $U$ be any vertical
vector field on $X$. Then we have $\pi_*[Y, U]=[\pi_*Y, \pi_*U]=0$, namely $[Y, U]$ is a 
vertical vector field on $X$. Thus one has, by the torsion-freeness of the connection $\nabla$,
\begin{eqnarray*}
&&<A(Y, Y), U> = <\nabla_YY, U> = -<Y, \nabla_YU>\\ 
&=& -<Y, [Y, U]+\nabla_UY> = -<Y, \nabla_UY> = -{\frac 1 2}U|Y|^2.
\end{eqnarray*}
Since $Y$ is the lift of a vector field $Y$ of $B$, so $|Y|^2$ is constant on any fiber 
of $\pi: X\to B$, and hence $<A(Y, Y), U>=0$ for any vertical vector $U$.
On the other hand, by the definition, $A(Y, Y)$ is a vertical vector, so $A(Y, Y)=0$. 

We now turn to the proof of the lemma. Take a normal frame $\{e_i\}$ of $M$ and a orthogonal frame of $X$ as follows: 
$\{e_{\alpha}, e_{\beta}, e_{\gamma}, \cdots,  e_{\mu}, e_{\nu}, \cdots\}$ 
with the properties $\{e_{\alpha}, \cdots\}$ being horizontal
and $\{e_{\mu}, \cdots\}$ vertical (note that, under such a restriction, one cannot get a normal frame 
in general). Then, under these frames, the stress-energy tensor of 
the horizontal map $u$ can be written as 
\[
\tau(u) = \sum_i\nabla{\text d}u(e_i, e_i) = \sum_{i, \alpha}u^{\alpha}_{ii}e_{\alpha} + 
\sum_{i, \alpha, \beta}u^{\alpha}_iu^{\beta}_i\nabla_{e_{\beta}}{e_{\alpha}}, 
\]
and hence its vertical part is 
$\sum_{i, \alpha, \beta}u^{\alpha}_iu^{\beta}_i{\cal V}\nabla_{e_{\beta}}{e_{\alpha}}$, 
which, by the previous discussion, is just
\[
{\cal V}\tau(u) = \sum_{i, \alpha, \beta}u_i^{\alpha}u_i^{\beta}A(e_{\alpha}, e_{\beta})
= \sum_{i}A({\text d}u(e_i), {\text d}u(e_i)) = 0.
\]
This completes the proof of the lemma.
 
The lemma 1 tells us that the solution $u(\cdot, t)$ to $(*)$ is horizontal for any $t\in [0, T)$, 
so ${\cal V}\tau(u(\cdot, t))=0$ for $t\in [0, T)$. Thus, by the previous computation, we actually obtain
\begin{eqnarray*}
(\Delta - {\frac{\partial}{\partial t}})e(u)
&=& |\nabla{\text d}u|^2 
+ <{\text{Ric}}^M({\text d}u(e_i), {\text d}u(e_i))>\\
&&- <R^X({\text d}u(e_i), {\text d}u(e_j)){\text d}u(e_i), {\text d}u(e_j)>.
\end{eqnarray*}
By the assumption on $\pi: X\to B$, $X$ has non-positive sectional curvature in the horizontal direction,
so we have
\[
(\Delta - {\frac{\partial}{\partial t}})e(u) \ge ce(u),
\]
for some constant $c$, which only depends on $M$. Denote the total energy of $u(\cdot, t)$
by $E(u(\cdot, t))$ for $t\in [0, T)$, i.e. $E(u(\cdot, t))=\int_Me(u(\cdot, t)){\text d}x$.
Then, one has
\begin{eqnarray*}
&&{\frac{d}{dt}}E(u(\cdot, t)) = {\frac{d}{dt}}\int_M<du, du>dx
= \int_M<\nabla_{\frac{\partial}{\partial t}}du, du>dx \\
&=& \int_M<\nabla{\frac{\partial u}{\partial t}}, du>dx
= -\int_M<{\frac{\partial}{\partial t}}u, \tau(u)>dx = -\int_M|{\cal H}\tau{(u)}|^2dx\le 0. 
\end{eqnarray*}
Summing all the above up, we have
\begin{lem}
Suppose $u(x, t)$ is a solution of $(*)$. Then
for some constant $c$, 
\[
(\Delta - {\frac{\partial}{\partial t}})e(u) \ge ce(u);
\]
furthermore, the total energy $E(u(\cdot, t))$ is a decreasing function of $t$.
\end{lem}

Combining the above lemma  with Lemma 2.3.1 in \cite{j}, one has
\begin{lem}
Let $t>0$, $0<R<\min(i(M), {\frac{\pi}{2\Lambda}})$, where $i(M)$ is the injective radius of $M$, 
and $\Lambda^2$ is an upper bound for the sectional curvature of
$M$. Then, for all $x\in M$, 
\[
e(u)(x, t)\le c (tR^{-m-2}+t^{-{\frac m 2}})\int_Me(g)(y){\text d}y,
\]
where $m={\text{dim}}M$ and $c$ is some constant depending only on the
geometry of $M$;  
and for any $t_0<t$, in particular $t_0=0$,
\[
e(u)(x, t)\le c R^{-2}\sup_{x\in M}e(u)(y, t_0).
\]
\end{lem}

In the following, we want to derive a stability lemma. Let $g(x, s)$ be a smooth horizontal family of
smooth horizontal maps from $M$ to $X$ with parameter $s\in [0, s_0]$, 
i.e. both $g(\cdot, s)$ for any $s\in [0, s_0]$ 
and ${\frac{\partial{g}}{\partial s}}$ being horizontal. 
Suppose that $u(x, t, s)$ is a family of solutions of $(*)$ with initial data $g(x, s)$ 
for $0\le s\le s_0$. As pointed out before, ${\frac{\partial{u}}{\partial t}}$ is horizontal;
using the same discussion as in Lemma 1, we now show
that ${\frac{\partial{u}}{\partial s}}$ is also horizontal: Fixing $x\in M$ and $s_1\in [0, s_0]$, one can then
consider $u(x, t, s)$ as a variation of the curve $u(x, t, s_1)$. Take
arbitrarily a vertical tangent vector $V$  
at $u(x, 0, s_1)$ and translate parallelly $V$ along the $t$-curve $u(x, t, s_1)$ 
and then the $s$-curves $u(x, t, s)$, still denoted by $V$. (Note that $V$ is not necessarily parallel 
along other $t$-curves $u(x, t, s)$ for $s\neq s_1$.) 
Compute ${\frac{\partial}{\partial t}}<{\frac{\partial u}{\partial
    s}}, V>$: 
\[
{\frac{\partial}{\partial t}}<{\frac{\partial u}{\partial s}}, V>=
<\nabla_{{\frac{\partial}{\partial t}}}{\frac{\partial u}{\partial s}}, V>
=<\nabla_{{\frac{\partial}{\partial s}}}{\frac{\partial u}{\partial t}}, V>=
{\frac{\partial}{\partial s}}<{\frac{\partial u}{\partial t}}, V>=0;
\]
on the other hand, $<{\frac{\partial u}{\partial s}}, V>|_{t=0}=0$, 
therefore $<{\frac{\partial u}{\partial s}}, V>=0$. 
Thus, for any fixed $t\in [0, T)$ and $s\in [0, s_0]$,
the derivative of $u$ with respect to $s$, ${\frac{\partial u}{\partial s}}(\cdot, t, s)$, can be 
considered as a horizontal vector field of $X$ (if necessary, it can be considered as
some section of a certain pull-back bundle). 
Using the horizontality of ${\frac{\partial u}{\partial s}}(\cdot, t, s)$, we then have
\begin{lem}
For every $s\in [0, s_0]$, the quantity
\[
\sup_{x\in M}|{\frac{\partial u}{\partial s}}|^2(x, t, s)
\]
is decreasing in $t$. Hence also the quantity
\[
\sup_{x\in M, s\in [0, s_0]}|{\frac{\partial u}{\partial s}}|^2(x, t, s)
\]
is a decreasing function in $t$.
\end{lem}
{\bf Proof.} As before, one can compute under a normal frame $\{e_i\}$
\begin{eqnarray*}
&&(\Delta-{\frac{\partial}{\partial t}})|{\frac{\partial u}{\partial s}}|^2\\
&=& 2|\nabla{\frac{\partial u}{\partial s}}|^2 - 
2\sum_i<R({\frac{\partial u}{\partial s}}, {\text d}u(e_i)){\frac{\partial u}{\partial s}}, {\text d}u(e_i)>.
\end{eqnarray*}
Here we use the heat eqaution and ${\frac{\partial u}{\partial s}}$'s horizontality. 
Thus, by the assumption on the sectional curvature in the horizontal
direction, we have 
\[
(\Delta-{\frac{\partial}{\partial t}})|{\frac{\partial u}{\partial s}}|^2\ge 0.
\] 
The lemma then follows from the maximum principle for parabolic
equations.

\vskip .3cm
In order to apply the regularity theorems for elliptic equations, we have to make sure that
the solution of $(*)$ with the given initial data $g$ lies in a suitable coordinate chart of 
$X$ when the domain considered is  small enough and the time interval
enough short. We have obtained  
a point-wise upper bound for the derivatives of $u$ with respect to
the space variables,  
so we still have to derive a bound for the time derivative of the solution. This can be done by 
applying the above lemma.
\begin{lem}
Suppose that $u(x, t)$ is a solution of $(*)$ with the initial data $g$ for $t\in [0, T)$. 
Then for all $t\in [0, T)$ and $x\in M$
\[
|{\frac{\partial u(x, t)}{\partial t}}|\le \sup_{y\in M}|{\frac{\partial u(y, 0)}{\partial t}}|.
\]
\end{lem} 
{\bf Proof.} Setting $u(x, t, s)=u(x, t+s)$, then $u(x, t, s)$ can be considered a family of solutions to 
$(*)$ with a family of initial data $u(x, s)$. Applying the preceeding lemma to $u(x, t, s)$, we then 
get the present lemma. 

\vskip .3cm
Fix $x\in M$ and $t\in [0, T)$. As before, we take a normal frame $\{e_i\}$ at $x$, and a orthogonal
frame $\{e_{\alpha}, e_{\beta}, e_{\gamma}, \cdots,  e_{\mu}, e_{\nu}, \cdots\}$ at $u(x, t)$ 
with the property that $\{e_{\alpha}, \cdots\}$ are horizontal
and $\{e_{\mu}, \cdots\}$ are vertical; as pointed out before, under such a restriction, 
one cannot get a normal frame at 
$u(x, t)$ in general. Then, the heat equation $(*)$ can be rewritten under such frames at $(x, t)$ as
\[
(*')~~~~~~~~~~~~~~~~~~~~~~~~~~~~~\sum_iu^{\alpha}_{ii} + 
\sum_{i,\beta,\gamma}\Gamma_{\beta\gamma}^{\alpha}u^{\beta}_iu^{\gamma}_i 
= {\frac{\partial u^{\alpha}}{\partial t}}.~~~~~~~~~~~~~~~~~~~~~~~~~~~~~~~~~~~~~~~
\] 
{\bf Remark.} Note that the solution with the initial data $g$ is horizontal 
for both the space variable and the time variable,
as seen in Lemma 1. So by Lemma 2, ${\cal V}\tau(u)=0$, i.e. ${\cal H}\tau(u)=\tau(u)$. 
Thus we can actually omit $\cal H$ in the equation $(*)$ and think that $u$ just satisfies 
the usual heat equation for harmonic maps, $\tau(u)-{\frac{\partial
    u}{\partial t}}=0$. 
In the following estimate, we will actually adopt this point of view 
although it will not be pointed out explicitly. 

\begin{lem}
Suppose that $u(x, t)$ is a solution of $(*)$ (or $(*')$) with the initial data $g$ for $t\in [0, T)$. 
Then for every $\alpha\in(0, 1)$
\[
\Vert u(\cdot, t)\Vert_{C^{2+\alpha}(M; X)} + 
\Vert{\frac{\partial u}{\partial t}}(\cdot, t)\Vert_{C^{\alpha}(M; X)} \le c,
\]
where $c$ depends on $\alpha$, the initial data $g(x)$, and the geometry of $M$ and $X$, 
but not on $t$.
\end{lem}
{\bf Proof.} Rewrite $(*')$ as 
\[
\sum_iu^{\alpha}_{ii} =
-\sum_{i,\beta,\gamma}\Gamma_{\beta\gamma}^{\alpha}u^{\beta}_iu^{\gamma}_i 
+{\frac{\partial u^{\alpha}}{\partial t}}.
\]
If we restrict the solution $u$ to a suitable small coordinate chart at the point 
$x_0\in M$, say $B(x_0, \rho)$ with $\rho$ enough small, and a suitable small time 
interval $[t_0, t_1]$, $u(x, t)$ will stay in a certain coordinate
chart of $X$  
by the lemma 4 and the lemma 6; moreover, those two lemmata also imply
that the right-hand 
side of the above equation is bounded (note that the bound does not depend on $t$), 
this, by the elliptic regularity theory, then implies a bound (again not depend on $t$) 
for $\Vert u(\cdot, t)\Vert_{C^{1+\alpha}(M; X)}$ on a smaller coordinate chart,
say $B(x_0, {\frac{\rho}{2}})$ (see \cite{j}, Theorem 2.2.1). 
Thus, the right-hand side of the following parabolic equation
\[
{\frac{\partial u^{\alpha}}{\partial t}} - \sum_iu^{\alpha}_{ii} =
\sum_{i,\beta,\gamma}\Gamma_{\beta\gamma}^{\alpha}u^{\beta}_iu^{\gamma}_i 
\]
is bounded (the bound being independent of $t$) in $C^{\alpha}(M; X)$, and hence the 
Schauder estimate for parabolic equations then implies the estimate in the lemma, 
at least in the above small coordinate chart; but $M$ is compact, so the estimate is valid
on $M$.

\vskip .3cm
Based on the local existence for solutions and the above Schauder estimate, 
one has the following global existence theorem for $(*)$ with the initial data $g$. 
\begin{thm}
The solution $u(x, t)$ of the heat equation $(*)$ with the horizontal initial data $g$ 
exists for all $t\in [0, \infty)$,
if the Riemannian submersion $\pi: X\to B$ satisfies the Chow condition and 
 $B$ has non-positive sectional curvature. 
\end{thm}

In the following, we will show that 
the global solution $u(\cdot, t)$ in the theorem above converges to 
a horizontal harmonic map as $t$ goes to infinity. 
As seen before, we have shown the energy decay formula, namely
\[
{\frac{\text d}{{\text d}t}}E(u(\cdot, t)) 
= -\int_M|{\frac{\partial u(x, t)}{\partial t}}|^2{\text d}x = -\int_M|{\cal H}\tau{(u)}|^2{\text d}x; 
\]
observe also that the energy function $E(u(\cdot, t))$ in $t$ is nonnegative 
for $t\in [0, \infty)$, so there exists a sequence $\{t_n\}_{n=1}^{\infty}$ with $t_n\to\infty$ 
as $n\to\infty$ satisfying ${\frac{\text d}{{\text d}t}}E(u(\cdot, t))|_{t_n}\to 0$ as $n\to \infty$, 
this is just equivalent to $\int_M|{\frac{\partial u}{\partial t}}(x, t_n)|^2{\text d}x\to 0$ as $n\to\infty$. 
On the other hand, as seen in Lemma 7, ${\frac{\partial u}{\partial t}}(\cdot, t)$ 
has a $C^{\alpha}$-bound independent of the time $t$, so we obtain
\begin{lem}
There exists a sequence $\{t_n\}_{n=1}^{\infty}$ with $t_n\to \infty$ as $n\to\infty$, for which
${\frac{\partial u}{\partial t}}(x, t_n)$ converges to zero uniformly in $x\in M$ as $n\to\infty$.
\end{lem}  

Lemma 7 also tells us that $u(\cdot, t)$ has a time-independent $C^{2+\alpha}$-bound, so one obtains, 
by possibly passing to a subsequence of $\{t_n\}$, that $u(\cdot, t_n)$ converges 
at least $C^2$-uniformly to a map $u: M\to X$, which then is also
horizontal; furthermore, 
since $\{u(\cdot, t_n)\}$ is at least $C^2$-uniformly convergent to $u$ and both 
$u(\cdot, t_n)$ and $u$ are horizontal, so by the Hopf-Rinow theorem, as mentioned in the Introduction, 
some $u(\cdot, t_n)$, and hence $g(\cdot)$, is homotopic to $u(\cdot)$
by some  
horizontal homotopy $h(\cdot, s)$ for $s\in [0, 1]$ with $h(\cdot, 0)=u(\cdot, t_n)$ 
and $h(\cdot, 1)=u(\cdot)$.
Here by the homotopy $h(\cdot, s)$ being horizontal 
we mean that $h(\cdot, s)$ for each $s\in[0, 1]$ is a horizontal map
and the $s$-curves are  
also horizontal. Again since $\{u(\cdot, t_n)\}$
uniformly converges to $u$, w.l.o.g., we can assume that the lengths of the $s$-curves $h(x, s)$ have 
a sufficiently small upper bound $\epsilon>0$ independent of $x\in M$.    
Now, consider the family of the solutions $u'(x, t, s)$ ($s\in [0, 1]$)
to $(*)$ with $h(x, s)$ as the family of initial maps. It is clear that $u'(x, t, 1)=u(x)$ 
since ${\cal H}\tau(u)=0$ and $h(x, 1)=u(x)$; while $u'(x, t, 0)=u(x, t+t_n)$. By the Lemma 5,
the supermum with respect to $x$ of the length of $s$-curves $u'(\cdot, t, s)$ is a decreasing 
function in $t$ and hence less than $\epsilon$. Since $\epsilon$ is arbitrary, we have that $u(x, t)$ converges 
uniformly to $u(x)$ in $t$ in the sense of $C^0$, not only for a subsequence $\{t_n\}$.
Applying this to the heat equation $(*)$, one obtains
\[
{\cal H}\tau(u)=0.
\] 
Finally, the horizontality of $u(x)$ and the Lemma 2 tell us that ${\cal V}\tau(u)=0$, and hence
\[
\tau(u) = 0,
\]
i.e. the limit $u$ is a horizontal harmonic map. Thus we have
\begin{thm}
Suppose that $\pi:X\to B$ is a Riemannian submersion satisfying the Chow conditions and
that $B$ has non-positive sectional curvature. Let $M$ be a compact Riemannian manifold and
$g: M\to X$  a horizontal smooth map from $M$ to $X$. Then there exists a 
horizontal harmonic map $u: M\to X$ from $M$ into $X$ that is
homotopic to $g$ by a  
horizontal homotopy.
\end{thm}

\noindent
{\bf Remark.} The theorem above is actually valid in a more general setting, 
namely the equivariant one: Let $\phi: \pi_1(M)\to \pi_1(X)$ be a homomorphism and 
$g$  a $\phi$-equivariant map from $M$ into $X$, then one can solve the corresponding heat equation $(*)$ and obtain similar results, e.g. the existence for $\phi$-equivariant horizontal harmonic maps. We omit this, 
but point out that in applications we shall just use that setting. We will come back to this in \cite{jy}.

\bigskip
\noindent
J\"urgen Jost:\\
Max-Planck-Institute for Mathematics in the Sciences, Leipzig, Germany\\
and\\
\noindent
Yi-Hu Yang:\\
Department of Applied Mathematics, Tongji University, Shanghai, China\\
{\it e-mail}: yhyang@mail.tongji.edu.cn
\end{document}